\documentclass[oneside,a4paper,12pt]{article}
\usepackage{mathtools}
\usepackage{amsmath}
\allowdisplaybreaks[4]
\usepackage{csquotes}
\usepackage{amsthm}
\usepackage{amssymb}
\usepackage{mathrsfs}
\usepackage{color}
\usepackage{bm}
\usepackage{comment}
\usepackage{fancyhdr}
\usepackage{geometry}
\usepackage{hyperref}
\usepackage{enumerate}
\usepackage{graphicx}
\usepackage{subfig}
\usepackage{float}
\usepackage{bbm}

 \usepackage{tikz}

\theoremstyle{definition}
\newtheorem{thm}{Theorem}[section]
\newtheorem{theorem}{Theorem}[section]
\newtheorem{corollary}[thm]{Corollary}
\newtheorem{lemma}[thm]{Lemma}
\newtheorem{claim}[thm]{Claim}
\newtheorem{proposition}[thm]{Proposition}
\newtheorem{definition}[thm]{Definition}

\newtheorem{example}[thm]{Example}
\numberwithin{equation}{section}

\newcommand{\R}{\mathbb R}

\renewcommand{\H}{\mathbb{H}}
\renewcommand{\S}{\mathbb{S}}

\date{}
\usepackage{fancyhdr}
\begin{document}
	
	\pagestyle{fancy}                   
	\fancyhead[RO,LE]{}   
	\chead{    }      
	\fancyhead[LO,RE]{ }

	\renewcommand{\headrulewidth}{0mm} 
	\renewcommand{\footrulewidth}{0mm} 


\title{\textbf{The polyhedral decomposition of cusped hyperbolic $n$-manifolds with totally geodesic boundary}}

\author{\medskip Huabin Ge, Longsong Jia, Faze Zhang}

\date{}
\maketitle
\begin{abstract}

Let $M$ be a volume finite non-compact complete hyperbolic $n$-manifold with totally geodesic boundary. We show that there exists a polyhedral decomposition of $M$ such that each
cell is either an ideal polyhedron or a partially truncated polyhedron with exactly one truncated face.  This result parallels Epstein-Penner's ideal decomposition \cite{EP} for cusped hyperbolic manifolds and Kojima's truncated polyhedron decomposition \cite{Kojima} for compact hyperbolic manifolds with totally geodesic boundary. We take two different approaches to demonstrate the main result in this paper. We also show that the number of polyhedral decompositions of $M$ is finite.
\end{abstract}

\maketitle




\section{Introduction}
In general, hyperbolic manifolds refer to Riemannian manifolds with constant sectional curvature $-1$. Thurston's geometrization conjecture states that compact 3-manifolds can be decomposed into pieces with geometric structures, most of which is hyperbolic. In order to study the geometry and topology of 3-manifolds, a useful way is to decompose a 3-manifold into basic blocks. For different purposes, common basic blocks are handlebody, polyhedra, etc.
We concern the polyhedra decomposition in this paper.

On the one hand, hyperbolic ideal tetrahedra can be pasted isometrically to obtain hyperbolic 3-manifolds, this method was creatively developed by Thurston and became the main method for constructing cusped 3-manifolds. On the other hand, it is natural to ask whether each cusped 3-manifold can be broken down into ideal tetrahedra. This is still an open question, often referred to as the Thurston geometric ideal triangulation conjecture, because Thurston once proved the famous hyperbolic Dehn filling theorem on the premise that it was correct. To some extent, it stands as one of the foremost unresolved problems in the field of three dimensional geometry and topology since the Thurston's Geometrization Conjecture and the Virtual Haken Conjecture were solved.

The geometric ideal triangulation conjecture, pointed by Gu\'{e}ritaud-Schleimer \cite{GS}, ``is a difficult problem in general. General results are known only when $M$ is restricted to belong to certain classes of manifolds: punctured-torus bundles, two-bridge link complements, certain arborescent link complements and related objects, or covers of any of these spaces". Some of the works we have learned about include \cite{Aki1}-\cite{Aki-SWY-2}, \cite{Gu}-\cite{Gue-2}, \cite{Ham-P}, \cite{Jorgen}, \cite{Lackenby}, \cite{Nimer} (of course, there may be some important related work, and we apologize for not noticing it). Moreover, it is shown that the geometric ideal triangulation conjecture is true virtually, i.e. there is a finite cover of the manifold which admits a geometric ideal triangulation~\cite{LST,FHH}. For compact three-dimensional manifolds with boundary, under some combinatorial condition, the first author and collaborators \cite{Ge1} show that the geometric ideal triangulation conjecture is true by using combinatorial Ricci flow methods.

Although the geometric triangulation conjecture is still widely open, there have been great breakthroughs about geometric polyhedral decompositions for volume finite hyperbolic manifolds. Epstein-Penner \cite{EP} constructed an ideal polyhedral decomposition for volume finite, non-compact complete hyperbolic $n$-dimensional manifolds. Later, Kojima established a variant of the Epstein-Penner decomposition, that is, decomposing a compact hyperbolic manifold with nonempty geodesic boundary into truncated polyhera, and gave a more visual construction of the decomposition for the special case of compact hyperbolic 3-manifolds with totally geodesic boundary \cite{Kojima}.
Subsequently, a very natural question arises: is there a similar decomposition for volume finite, non-compact hyperbolic manifolds with both cusps and totally geodesic boundaries? In this article, we give an affirmative answer, that is:



\begin{theorem}
\label{main}
Let $M$ be a volume finite, non-compact, complete hyperbolic $n$-manifold with both cusps and totally geodesic boundaries. Then $M$ admits an ideal polyhedral decomposition with each cell either an ideal polyhedron or a partially truncated polyhedron with exactly one  truncated face.
\end{theorem}

We take two different approaches to demonstrate the above result. The first proof is a modification of \cite{EP}: The hyperbolic double $DM$ obtained by doubling $M$ is a volume finite, non-compact complete hyperbolic $n$-dimensional manifold with at least one cusp. By Epstein-Penner \cite{EP} theory, an ideal polyhedral decomposition of $DM$ is obtained. Then by carefully examining the symmetry properties of the decomposition along $\partial M$, we finally obtain a polyhedral decomposition of $M$. The second proof adopts a similar strategy of \cite{Kojima}, but is more straightforward.  We first construct a cell decomposition of the cut locus and then construct a geometric decomposition of $DM$ dual to the decomposition of cut locus. Finally, by the symmetry properties of the decomposition along $\partial M$, we obtain a polyhedral decomposition of $M$. One noticeable difference is that the first method deals with any dimension $n$, and in order to get a more visual decomposition, the second method focuses on dimension 3.

The two approaches have their own advantages and complement each other, the first method facilitates the direct handling of arbitrary dimensions, and the second method provides better operational visibility. For the dimension $n=3$, the two methods are essentially equivalent, and they are realized by examining the hyperbolic double $DM$ of $M$ along its geodesic boundaries $\partial M$ under different hyperbolic space models. Furthermore, by appropriately modifying the second method, a proof for any dimension $n$ can also be obtained. 

Based on Epstein-Penner \cite{EP}'s construction of polyhedral decompositions on cusped hyperbolic $n$-manifolds, Akiyoshi \cite{Aki2} showed that if a cusped $n$-manifold $N$ has at least two cusps, then by choosing decoration at each cusp with different volumes, the number of ideal polyhedral decompositions on $N$ is finite. In this paper, since $M$ is a volume-finite hyperbolic $n$-manifold that has both cusps and boundaries, $DM$ is a volume-finite cusped hyperbolic $n$-manifolds with at least two cusps. Therefore, by Akiyoshi \cite{Aki2}'s results and the first proof of Theorem \ref{main}, we have the following corollary:

\begin{corollary}
Let $M$ be a volume finite, non-compact, complete hyperbolic $n$-manifold with both cusps and totally geodesic boundaries. Then the number of ideal polyhedral decompositions of $M$ is finite.
\end{corollary}


The paper is organized as follows. We first introduce some basic notions in Section \ref{2}, including various models of hyperbolic spaces, the ideal polyhedron and partially truncated polyhedron. In Section \ref{section:approach1}, we provide a method which is a modification of \cite{EP} to approach Theorem \ref{main}.
In Section \ref{Section:proof2}, we adopt Kojima's proof philosophy in \cite{Kojima} to approach Theorem \ref{main}. Finally, we compare the two proofs and look forward to some potential applications in Section \ref{section:comment}.

~

\noindent
\textbf{Acknowledgements:}
The authors are very grateful to Professor Ruifeng Qiu, Feng Luo, Tian Yang for many discussions on related problems in this paper. The first two authors would like to thank Professor Gang Tian for his constant encouragement and support. Huabin Ge is supported by NSFC, no.12341102, no.12122119. Faze Zhang is supported by NSFC, no.12471065.

\section{Preliminaries}\label{2}


\subsection{Hyperbolic space}
To get a clearer understanding of our construction process, this section will introduce several models of $n$-dimensional hyperbolic space $\H^n$. Certainly, these models of the hyperbolic space are all equivalent and well-known, see Benedetti-Petronio \cite{BenPet}, Marden \cite{Marden}, Martelli \cite{Martelli}, Purcell \cite{Purcell}, Ratcliffe \cite{Ratcliffe}, Thurston \cite{Thurston} for instance. The main purpose of this section is to unify the  notations.

(1) \emph{The hyperboloid model for} $\mathbb{H}^n$. In $\mathbb{R}^{n+1}$, we consider the standard symmetric bi-linear form of signature $(n,1)$, which is called the Lorentzian scalar product,
\begin{equation*}
\langle x, y\rangle=-x_{0}y_{0}+x_{1}y_{1}+\cdots+x_{n}y_{n},
\end{equation*}
where $x=(x_{0},x_{1},\cdots,x_{n})$ and $y=(y_{0}, y_{1}, \cdots, y_{n})$. $\mathbb{R}^{n+1}$ with such a Lorentzian scalar product is called the Minkowski space, and is often denoted by $\mathbb{R}^{n,1}$ in the literature. The following set
\begin{equation*}
V=\big\{x\in \mathbb{R}^{n,1}|-x_{0}^{2}+x_{1}^{2}+\cdots+x_{n}^{2}=-1\big\}
\end{equation*}
is a hyperboloid with two sheets. Denote $V^{+}$ by the connected component of $V$ with $x_{0}>0$. $V^{+}$ is called the hyperboloid model for hyperbolic space $\mathbb{H}^n$. The restriction of the Lorentzian scalar product to the tangent space $T_xV^+$ at each $x\in V^+$ is positive definite and hence induces a metric tensor on $V^+$. This makes $V^+$ a $n$-dimensional Riemannian manifold, which is simply connected, complete, and has a constant sectional curvature $-1$. Let $O(n, 1)$ be the group of all linear isomorphisms of $\mathbb{R}^{n, 1}$ that preserve the Lorentzian scalar product. An element in $O(n, 1)$ always preserves the two sheeted hyperboloid $V$, and the elements preserving the upper sheet $V^+$ form a subgroup of index two in $O(n, 1)$, which is denoted by $O^+(n, 1)$ and is isomorphic to the isometric group $\text{Isom}(V^+)$ of $V^+$.

(2) \emph{The unit ball model for} $\mathbb{H}^n$. This model is also called the Poincar\'{e} unit ball model, or the Poincar\'{e} model. Consider the unit ball in the $n$-dimensional Euclidean space $\mathbb {R}^{n}$ is
\begin{equation*}
B^{n}=\big\{x=(x_{1},\cdots,x_{n})\in \mathbb{R}^{n}\,\big|\,|x|<1\big\}.
\end{equation*}
The unit ball model for $\mathbb{H}^n$ is $B^{n}$ equipped with the following conformal metric:
\begin{equation*}
ds^{2}=\frac{4|dx|^{2}}{(1-|x|^{2})^{2}} .
\end{equation*}

(3) \emph{The upper half-space model for} $\mathbb{H}^n$. The upper half-space of $\mathbb {R}^{n}$ is
\begin{equation*}
H^{n}=\big\{x=(x_{1},\cdots,x_{n})\in \mathbb {R}^{n}\;\big|\,\,x_{n}>0\big\}.
\end{equation*}
The upper half-space model for $\mathbb{H}^n$ is $H^{n}$ equipped with the following metric:
\begin{equation*}
ds^{2}=\frac{|dx|^{2}}{x_{n}^{2}} .
\end{equation*}

(4) \emph{The projective model for} $\mathbb{H}^n$. Consider the real projective space $RP^{n}$. It is well-known that there is a canonical projection of $\mathbb{R}^{n,1}$ onto $RP^{n}$ \cite{Marden}\cite{Ratcliffe}\cite{Thurston}. The open unit ball $B^{n}$ of $RP^{n}$ quipped with the pull-back metric is called the projective model or the Klein model for $\mathbb{H}^n$. In this model, $\mathbb{H}^n$ is identified to the open unit ball in $\mathbb{R}^{n}\subset RP^{n}$, that geodesics of $\mathbb{H}^n$ then correspond to the intersection of straight lines of $\mathbb{R}^n$ with $\mathbb{H}^n$, and that totally geodesic planes in $\mathbb{H}^n$ are the intersection of linear planes with $\mathbb{H}^n$. Equivalent, a more intuitive explanation is as follows. Identifying the Poincar\'{e} unit ball model $(B^n, \frac{4|dx|^{2}}{(1-|x|^{2})^{2}})$ with the subset $\{x\in\mathbb{R}^{n+1}:x_1^2+\cdots+x_n^2<1, x_0=0\}$ in $\mathbb{R}^{n+1}$ with the corresponding metric, and mapping $B^n$ under the stereographic projection $\Pi$ with respect to the south pole $(-1, 0, \cdots, 0)$, we obtain the upper half of the unit sphere
$$S^n_{+}=\{x\in \mathbb{R}^{n+1}:x_0^2+x_1^2+\cdots+x_n^2=1, x_0>0\}.$$
Thus, composing $\Pi$ with the projection $P: (x_0, x_1, \cdots, x_n)\mapsto (x_1, \cdots, x_n)$ we obtain a homeomorphism $P\Pi$ from $B^n$ onto itself. $P\Pi$ extends continuously to the boundary $S^{n-1}$ of $B^n$ by setting $P\Pi(x)=x$. A geodesic line in the Poincar\'{e} model $(B^n, \frac{4|dx|^{2}}{(1-|x|^{2})^{2}})$ of $\mathbb{H}^n$ is mapped under $P\Pi$ to Euclidean segment with the same end poingts. The unite ball $B^n$ with the metric induced by $P\Pi$ from $(B^n, \frac{4|dx|^{2}}{(1-|x|^{2})^{2}})$ is the Klein model of $\mathbb{H}^n$. For instance, see the following Figure \ref{Fig-Klein} for $n=2$ case.
\begin{figure}[htbp]
\centering
\includegraphics[scale=0.4]{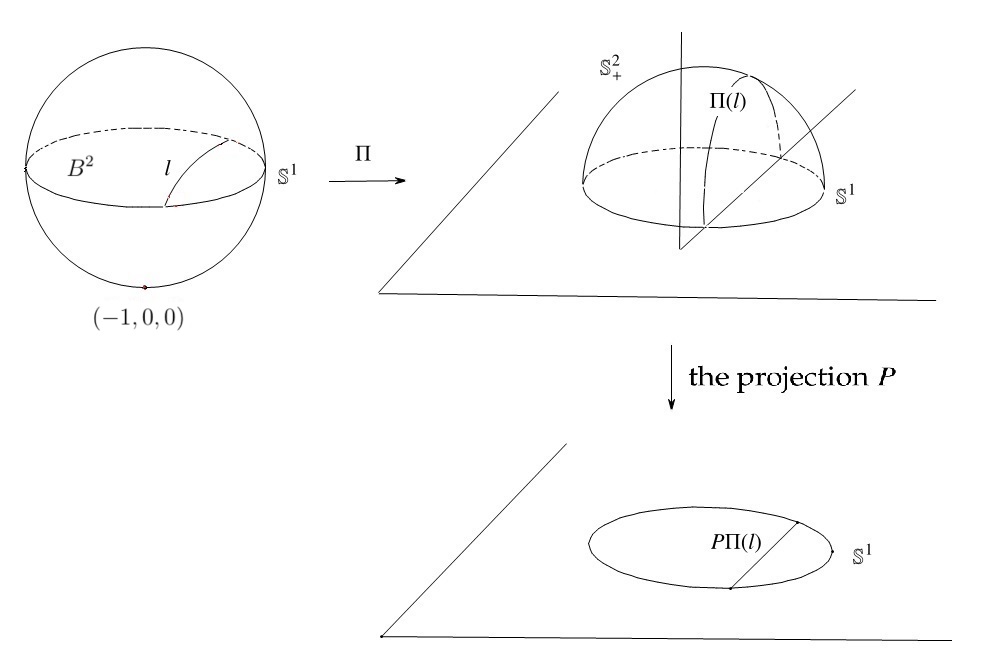}
\caption{The stereographic projection and the Klein model $B^2$}
\label{Fig-Klein}
\end{figure}

It is well-known that the above four different models for $\mathbb{H}^n$ are isometrically diffeomorphic to each other. In addition, both the unit ball model and the projective model are represented by $B^n$. Since we will deliberately point out the meaning of $B^n$ whenever it appears later, the reader will not be confused.

\subsection{ Partially truncated polyhedron}\label{subsection:2.2}

In this subsection, we use the Klein model $B^{n}$ of $\mathbb{H}^n$ and assume that $B^{n}\subset RP^{n}$. By Luo-Schleimer-Tillmann~\cite{LST}, we have the following definitions.
\begin{definition}
 A \textbf{projective polyhedron of $p$-$q$ type} $\hat{P}$  is defined as a convex projective polyhedron with $p+q$ vertices of $RP^{n}$ such that
  \begin{itemize}
    \item[(1)]  $p$ vertices of $\hat{P}$ lie in $RP^{n}\setminus \overline{B^{n}}$ and the other $q$ vertices of $\hat{P}$ lie in $\partial B^{n}$;
    \item[(2)]   Each $n-2$ dimensional face of $\hat{P}$ intersects $B^{n}$ non-empty. The $p$ vertices denoted by $v$ lying in $RP^{n}\setminus \overline{B^{n}}$ and the $q$ vertices lying in $\partial B^{n}$ are called \textbf{hyperideal vertex} and \textbf{ideal vertices} of $\hat{P}$, respectively.
 \end{itemize}

\end{definition}

For any hyperideal vertex $v$ of a projective polyhedron of $p$-$q$ type $\hat{P}$. There is associated hyperplane $H(v)$ such that
 \begin{itemize}
    \item[(1)] $H(v)$ is parallel to the Euclidean orthogonal complement of $v$;
    \item[(2)] $H(v)$ meets $\partial B^n$ in the set of all points $x$ such that there is a tangent line to $\partial B^n$ passing through $x$ and $v$.
 \end{itemize}

\begin{definition}

    \textbf{Truncation} of a projective polyhedron of $p$-$q$ type $\hat{P}$ at a hyperideal vertex $v$ is defined as cutting off an open star of $v$ from $\hat{P}$ along hyperplane $H(v)$. Convex domain $P$ obtained by removing the $q$ ideal vertices and truncating $q$ hyperideal vertices is called a \textbf{partially truncated polyhedron of $p$-$q$ type}. The boundary faces of $P$ resulting from truncation are called \textbf{external face} and the other boundary faces are called \textbf{internal faces}.
\end{definition}

\begin{figure}[htbp]
\centering
\includegraphics[scale=0.27]{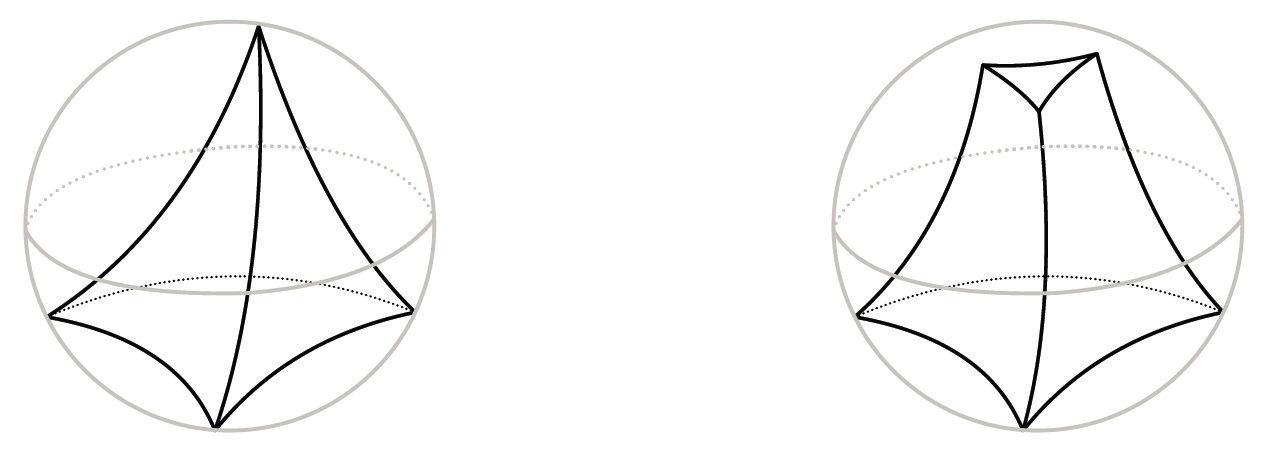}
\caption{An ideal tetrahedron / a type 1-3 partially truncated tetrahedron in $\mathbb{H}^3$}
\label{Fig-tetra}
\end{figure}

When $p=0$, a partially truncated polyhedron of $p$-$q$ type is often called an \textbf{ideal polyhedron}. We remark that, for dimension three, by Rivin~\cite{Rivin} and Bao-Bonahon \cite{BB}, given the combinatorial structure, the shape of an $3$-dimensional (hyper-)ideal polyhedron $P$ is fully determined by its dihedral angles. For instance, see Figure \ref{Fig-tetra} for an ideal tetrahedron and a type 1-3 partially truncated tetrahedron in $\mathbb {H}^{3}$.

Note for a partially truncated polyhedron, the external face intersects each internal faces orthogonally with respect to the pull-back metric in $RP^n$. 

\begin{definition}
Let $M$ be a volume finite, non-compact hyperbolic $n$-dimensional manifold with totally geodesic boundary. A \textbf{mixed ideal polyhedra decomposition} of $M$ refers to a decomposition of $M$ where each element is either an ideal polyhedron or a partially truncated polyhedron of $p$-$q$ type. If all elements is ideal polyhedron, i.e. partially truncated polyhedron of $0$-$q$ type, the decomposition is called an \textbf{ideal polyhedra decomposition}.
\end{definition}


\section{The First Proof to Theorem  \ref{main}}
\label{section:approach1}

In this section, we present the first method for proving Theorem \ref{main}. The proof involves considering the hyperbolic double $DM$ and applying the construction of Epstein-Penner \cite{EP}. The key point is to analyze the symmetry properties that arise in the decomposition of $DM$, which is precisely stated in Propositions \ref{lemma1}, \ref{lemma2}, \ref{lemma:symmetricLine} and \ref{lemma3}.

The construction process consists of three steps, and is mainly described in the hyperboloid model and supplemented by the Klein model. The first step is to choose decorations on cusps of $DM$ and identify their preimages in the universal covering as points on the light cone. Next step is to take the Euclidean convex hull of these points. The final step is to project the convex hull vertically onto the hyperboloid to obtain a $\pi_1(DM)$-invariant decomposition which can be passed onto the hyperbolic double $DM$.  The geodesic boundary of $M$ are considered as embedded surfaces of $DM$. Along the above process, we need to show some symmetry properties for the decomposition of $DM$ along the totally geodesic boundaries. These symmetry properties allow us to obtain a decomposition of $M$ itself.

This section will be divided into three parts, corresponding to the above three steps of the construction.

\subsection{Identification of horoballs on the light cone}
Let $M$ be a volume finite, non-compact, complete hyperbolic $n$-manifold with totally geodesic boundary. Assume $M$ has $k$ boundary components, which is denoted by $W_{1}$, $W_2$, $\cdots$ and $W_k$. We construct a new manifold $DM$ by pasting $M$ together with its mirror image along all boundary components of $M$. The resulting hyperbolic $n$-manifold $DM$ is called the \emph{hyperbolic double} of $M$  with respect to the whole boundary of $M$.

$DM$ is a volume finite, non-compact hyperbolic $n$-manifold without boundary. Moreover, the geodesic boundary components of $M$ are closed totally geodesic submanifolds of $DM$, which are still denoted by $W_{1}$, $W_2$, $\cdots$, $W_k$. Notably, $DM$ is symmetric with respect to $\{W_{j}\}$, that is, each point $x\in M$ and its mirror image $x'$ are symmetric with respect to $\{W_{j}\}$ along the minimal geodesic lines toward each $W_j$, $1\leq j\leq k$.

Since $DM$ is a complete and hyperbolic $n$-manifold, it has an universal covering space $\mathbb{H}^n$ with $\mathbb{H}^n/\Gamma=DM$, where $\Gamma\cong\pi_1(DM)$ is the covering transformation group and is discrete. Let $\pi:\H^n\rightarrow DM $ be the universal covering map. Since $DM$ is symmetric with respect to each $\{ W_j\}$ for $1\leq j\leq k$, we obtain a similar symmetry property for the universal covering:

\begin{lemma}
\label{lemma1}
Let $W_j$ be a connected component of the totally geodesic boundary of $M$ and $W^l_{j}$ be a component of the pre-image of $W_{j}$ under the universal covering map. Denote the symmetry transformation in $\H^n$ with respect to $W^l_{j}$ by $\tau$. Then $\tau^{-1}\Gamma \tau=\Gamma$.
\end{lemma}

\begin{proof}
For each connected component $W_j$ of the totally geodesic boundary of $M$ and each component $W_j^l$ of its pre-image, since $DM$ is symmetric with respect to $W_{j}$, there exists an isometry $\tau'$ of $\H^n$, which is the lift of the hyperbolic symmetry $DM\rightarrow DM$ and satisfying $(\tau')^{-1}\Gamma (\tau')=\Gamma$. On the other hand, $\tau$ is the reflection in the hyperplane containing $W_j^l$, hence $\tau=\tau'$ on $W^l_{j}$. Moreover, for any point $q\in W_j^l$, $\tau$ equals to $\tau'$ on the fundamental domain containing $q$ and hence then on the tangent space of $q$. Thus $\tau=\tau'$ and the lemma is proved.
\end{proof}

Since $M$ is a non-compact and finite-volume, it contains at least one cusp. Then $DM$ contains at least two cusps and those cusp are symmetric along $\{W_j\}$ in pairs. By the symmetry property of $DM$ along $\{W_j\}$ and the Proposition $14.1$ of \cite{Purcell}, we can always choose a horoball neighborhood around each cusp of $DM$, so that
 \begin{itemize}
    \item[(1)] it is far away from $\partial M$;
    \item[(2)] if two cusps are symmetric along $
    \{W_j\}$, then the  horoball neighborhoods are also symmetric;
    \item[(3)] the liftings $\widetilde{B}$ in $\mathbb{H}^n$ of these horoball neighborhoods are disjoint horoballs.
 \end{itemize}
 The horoball neighborhood is refered to a \emph{decoration} of the corresponding cusp.

Recall that $V$ is the two-sheeted hyperboloid in $\mathbb{R}^{n,1}$. The \emph{light cone} $L$ of the hyperboloid $V$ is defined by
\begin{equation*}
L=\big\{x\in \mathbb{R}^{n,1}|x_{0}^{2}=x_1^2+\cdots+x_{n}^{2}\big\},
\end{equation*}
where $x=(x_{0},x_{1},\cdots,x_n)$. The positive light cone $L^{+}$ is the connected component of $L\setminus\{0\}$ satisfying $x_0>0$. Note that in $\mathbb{R}^{n,1}$, a ray from the origin in $L^{+}$ corresponds to a point on $\partial B^{n}$ where $B^{n}$ is the unit ball model of $\mathbb{H}^n$. Horoballs in $V^+$ can be described by points in $L^+$:
\begin{definition}[\cite{EP}]
A horoball in $V_+$ can be described as
\begin{equation*}
\{w\in \mathbb{H}^3|-1\leq \langle w,p\rangle <0\}.
\end{equation*}
Here $p\in L^{+}$ is referred as the $\textbf{center}$ of the horoball.
\end{definition}

Denote the centers of horoballs in $\widetilde{B}$ as $\widetilde{\mathcal{B}}$. Then $\widetilde{\mathcal{B}}$ can be decomposed into $\cup_{i=1}^{2\nu} \widetilde{\mathcal{B}_i}$ such that each $\widetilde{\mathcal{B}_i}$ according to the corresponding cusp of $DM$ where $\nu$ is the number of cusps in $M$. By Theorem $2.4$ of \cite{EP},  $\widetilde{\mathcal{B}}$ is discrete in $\mathbb{R}^{n,1}$.

By Lemma \ref{lemma1}, for any connected component of the preimage of a geodesic boundary, the symmetry map along it yields a symmetry map between $\widetilde{\mathcal{B}}$ which has the following property:

\begin{proposition}\label{lemma2}
Let $W_j$ be a connected component of the totally geodesic boundary of $M$, and $W^l_{j}$ be a component of the pre-image of $W_{j}$ under the universal covering map.
Then there exists a direction $ v \in \mathbb{R}^{n,1}$ that depends only on $W^l_{j}$ such that the line connecting the points $p_1$ and $p_2$ is parallel to $v$,  where $p_1, p_2\in \widetilde{\mathcal{B}}$ are arbitrary ideal points symmetric with respect to $W^l_{j}$.
\end{proposition}

\begin{proof}
We first consider a special case. If $W^l_{j}$ passes through the  $x_0$-axis, then $p_1$ and $p_2$ are symmetric along a hyperplane passing through  $x_0$-axis. Let $v$ be the vector perpendicular to the plane that passes through $W^l_{j}$ in Euclidean metric, then the line connecting the points  $p_1$ and $p_2$ in $L^{+}$ is parallel to $v$.

In general, if $W^l_{j}$ does not pass through the $x_0$-axis, there exists an isometric transformation  $A \in O^+(n,1)$ that depends only on $W^l_{j}$, such that $A(W^l_{j})$ passes through the $x_0$-axis.

Let
\begin{equation*}
B_1'= \{ w \in \mathbb{R}^{n} \mid \langle w, p_1 \rangle = -1 \}
\end{equation*}
\begin{equation*}
B_2'= \{ w \in \mathbb{R}^{n} \mid \langle w, p_2 \rangle = -1 \}
\end{equation*}
be the horoballs whose centers are $p_1$ and $p_2$ respectively.

By the isometric transformation  $A \in O^+(n,1)$, $B_1'$ is transformed to
\begin{equation*}
A(B_1') = \{ A(w) \mid \langle w, p_1 \rangle = -1 \}
\\= \{ w \in \mathbb{R}^{n} \mid \langle w, A^{-1}(p_1) \rangle = -1 \}.
\end{equation*}
So the center of the horoball $A(B_1')$ is $ A^{-1}(p_1) $.~By a similar argument for $B_2'$, the center of the horoball $A(B_2')$ is $ A^{-1}(p_2) $.~Thus,~by the argument in the special case,~the line connecting the points $A^{-1}(p_1)$ and $ A^{-1}(p_2)$ is parallel to a vector $v$,~and further the line connecting the points  $p_1$ and $p_2$ is parallel to the vector $ A(v) $.
\end{proof}

\subsection{Construction of the Euclidean convex hull}

Let $C$ be the closed convex hull in $\R^{n,1}$ of $\widetilde{\mathcal{B}}=\cup_{i=1}^{2\nu} \widetilde{\mathcal{B}_i}$. Then Epstein-Penner have the following proposition:

\begin{proposition} [\cite{EP}]
The boundary of $C$ in $\mathbb{R}^{n,1}$ is the union of $C\cap L^{+}$ and a countable set of $(n-1)$-dimensional faces $F_{1},F_{2},\ldots $, each of which is a convex hull of a finite number of points in $\widetilde{\mathcal{B}}$. Here, $C\cap L^{+}$ is the set of points of the form $\alpha z$,  where $\alpha\geq 1$ and $z\in \widetilde{\mathcal{B}_i}$ for $1\leq i \leq 2\nu$.
\end{proposition}

Moreover, we also have a symmetry property regarding the convex hull:~

\begin{proposition}\label{lemma:symmetricLine}
Let $W_j$ be a connected component of the totally geodesic boundary of $M$, and $W^l_{j}$ be a component of the pre-image of $W_{j}$ under the universal covering map $\pi$. Then the boundary of the convex hull as a $(n-1)$-dimensional Euclidean infinite polyhedron will be symmetric along $W_j^l$.
\end{proposition}
\begin{proof}
For any two points $p,q\in \widetilde{\mathcal{B}}$, let $p'$ and $q'$ be the symmetric points of $p$ and $q$, respectively,  along $W_j^l$.~By proposition \ref{lemma2}, the lines $pp'$ and $qq'$ are parallel.~Thus, the convex hull spanning by $p$, $q$, $p'$ and $q'$ will be symmetric along $W_j^l$. Since all the points of $\widetilde{\mathcal{B}}$ are symmetric in pairs, the boundary of the convex hull as a $(n-1)$-dimensional Euclidean infinite polyhedron will be symmetric along $W_j^l$.
\end{proof}

\subsection{Construction of the decomposition on $M$}

By projecting the convex hull onto the hyperboloid $V^+$ vertically, we will obtain a $\Gamma$-invariant ideal polyhedral decomposition of $\H^n$ which descends to an ideal polyhedral decomposition $\mathcal D$ of $DM$.

\begin{proposition}\label{lemma3}
The ideal polyhedral decomposition $\mathcal D$ of $DM$ obtained above will induce a mixed ideal polyhedra decomposition of $M$.
\end{proposition}
\begin{proof}
There are two cases to consider.

Case $1$: If there is an ideal polyhedron $P\in \mathcal D$ such that $P$ does not intersect with any embedded submanifold $W_{j}$ for $1\leq j\leq k$ in $DM$.

By Proposition \ref{lemma:symmetricLine}, there must exist another ideal polyhedron $P'$ that is symmetric to $P$ with respect to $W_{j}$. After the symmetrical action on $DM$ along $W_j$, we obtain the original manifold $M$, and the images of $P$ and $P'$ will be two ideal polyhedra of $M$.

Case $2$: If there is an ideal polyhedron $P''$ in $\mathcal D $ such that $P''$ intersects some embedded submanifold $W_{j}$ for $1\leq j\leq k$ in $DM$.

Note that the boundary of the convex hull $C$ gives a $\Gamma$-invariant tessellation of $\H^n$. Then by proposition \ref{lemma:symmetricLine}, for any component $W_j^l$ of the preimage of $W_j$, the image of the projection from $C$ onto the hyperboloid $V^+$ vertically will be symmetric with respect to $W_j^l$. Thus we obtain a tessellation of $\H^n$ where the regions, $n-2$ dimensional faces and vertices of the tessellation are all symmetric along $W_j^l$.

Through the covering map, $\H^n$ maps to $DM$. Consequently there is an ideal decomposition of $DM$ inherited from the tessellation. This ideal decomposition of $DM$
is symmetric along $W_j$ which means that the regions, $n-2$ dimensional faces and vertices of the the ideal decomposition are symmetric along $W_j$.

Furthermore, if a $n-2$ dimensional face in the ideal polyhedron of $DM$ intersects non-trivially with a geodesic surface $W_j $, then this face contains a fixed point under the symmetric map, indicating that this edge is symmetric with itself since symmetric map preserves the face. By that the symmetric map is a self-isometric transform of $DM$ and the face is geodesic, this face is perpendicular to $W_j$.

Hence, the symmetric map of  $DM$ along $W_j$ restricted to $P''$ is also a symmetric transform, and $W_j\cap P''$ is the set of all fixed points. Thus, $W_j\cap P''$ is a polygon that is perpendicular to all internal faces of $P''$.

By quotienting the symmetric action on $DM$, we obtain the quotient manifold $M$ and get a convex polyhedron $P'''$ in $M$. Moreover, in this case, the resulting polyhedron has only one boundary face such that whose all vertices are not ideal.

\begin{claim}\label{claim5}
This convex polyhedron $P'''$ is a partially truncated polyhedron of $1$-$m$ type in $\H^n$ where $m$ is a positive integer.
\end{claim}

\begin{proof}

Now place $P'''$ into the projective model $B^{n}\subset RP^{n}$ of $\H^n$, then these $n-2$ dimensional faces perpendicular to $W_j\cap P''$ is still perpendicular to $W_j$ with respect to the pull-back metric in $RP^n$. By \cite{Matthieu}, there is a point $v\in RP^{n}\setminus \overline{B^{n}}$ such that $H(v)=W_j$.

Then all $n-2$ dimensional subspaces containing $v$ are perpendicular to $W_j$ with respect to the pull-back metric in $RP^n$. Note that for any $n-2$ dimensional subspace $Z$ perpendicular to $W_j$,  $Z\cap W_j$ is a $n-3$ dimension subspace in $\R^n$. Then there is a $n-2$ dimensional subspace $Z'$ containing $Z\cap W_j$ and $v$. Both $Z$ and $Z'$ are perpendicular to $W_j$ at $Z\cap W_j$ with respect to the pull-back metric in $RP^n$ while there is exactly one $n-2$ subspace will satisfy the property. It follows that $Z=Z'$.

This implies that $W_j\cap P'''$ can be viewed as truncating $v$ in some projective polyhedron along $W_j$. And by that the other vertices of number $m$ which lie in $\partial H^{n}$ will be the ideal vertices of $P'''$,  we have $P'''$ is a partially truncated polyhedron of $1$-$m$ type for some positive integer $m$.
\end{proof}
Hence we finished the proof of the proposition \ref{lemma3}.
\end{proof}

From the above discussion, we finally obtain a mixed ideal polyhedra decomposition of $M$ such that each
cell is either an ideal polyhedron or a partially truncated polyhedron with exactly one hyperideal vertex. Thus Theorem \ref{main} is proved.

\begin{figure}[htbp]
    \centering
    \includegraphics[scale=0.84]{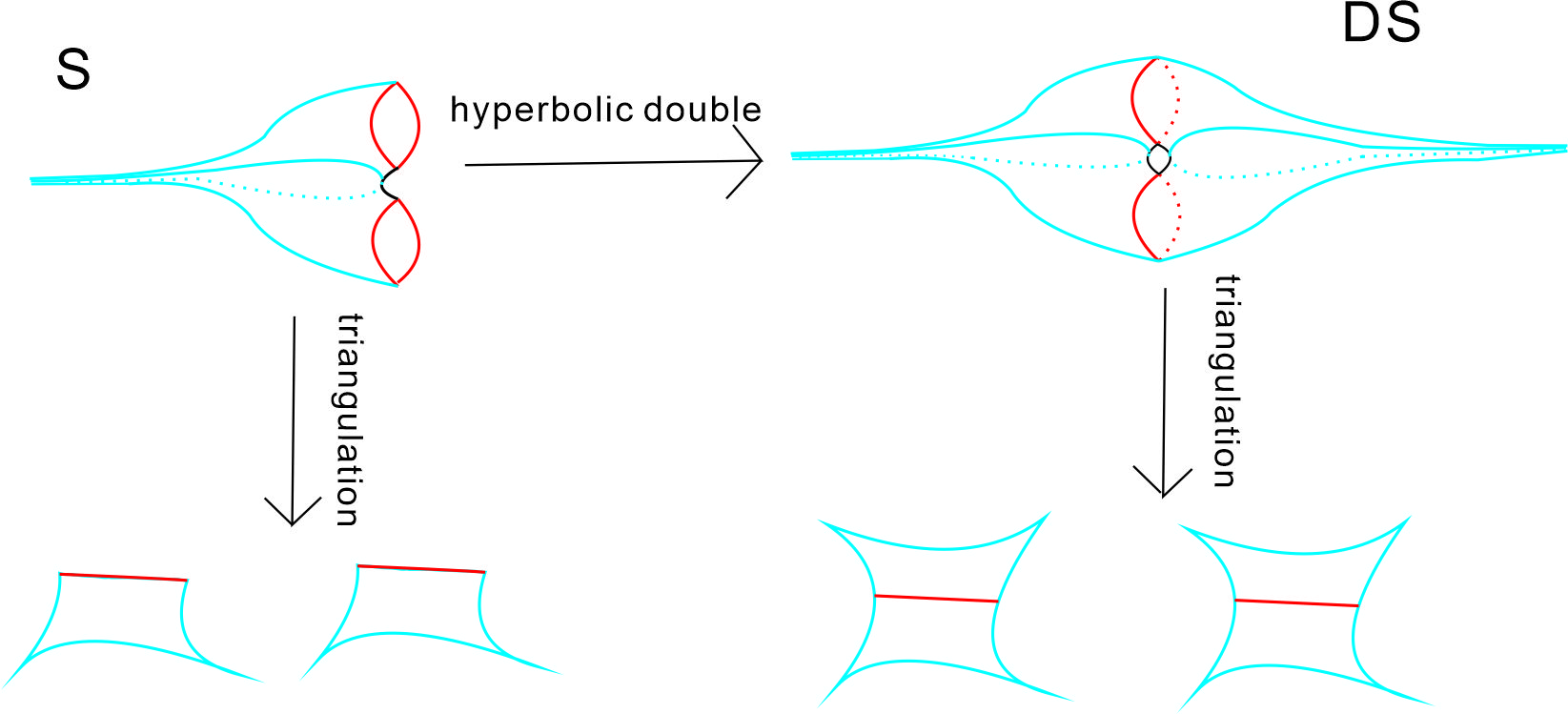}
    \caption{An example}
    \label{Fig3}
    \end{figure}

\begin{example}
We provide an example of $2$-dimensional manifold to illustrate our construction. Consider a one-cusped surface with two totally geodesic boundaries, which is denoted by $S$. Then the hyperbolic double along the boundaries is a two cusped surface without boundary, denoted by $DS$, see Figure \ref{Fig3}. The red circles represent the boundaries of $S$ and the blue lines indicate an ideal triangulation of $DS$. Then after applying the symmetric action, $S$ has a mixed ideal triangulation $\mathcal T$ such that each $2$-dimensional cell of $\mathcal T$ is a partially truncated polyhedron of $1$-$2$ type.
\end{example}


\section{The Second Proof to Theorem \ref{main}}\label{Section:proof2}

In this section, we will provide a proof of Theorem \ref{main} in the three-dimensional case. It should be noted that the resulting decomposition is essentially the same as the three-dimensional case in Section \ref{section:approach1}, but the construction is more geometric and intuitive.

The method used in this section are actually a modification of Kojima's construction of truncated polyhedral decompositions for compact $3$-manifolds with totally geodesic boundary. However, the difference between our method and Kojima's is that instead of studying the cut locus of the totally geodesic boundary, we will study the cut locus of decorations of cusps.

Let $M$ be a volume finite, non-compact, complete hyperbolic $3$-manifold with totally geodesic boundary and \{$W_{1}$, $W_2$, $\cdots$, $W_k$\} be all the boundary components of $M$. By pasting $M$ together with its mirror image along all boundary components of $M$, we get a new manifold $DM$ which is called the \emph{hyperbolic double} of $M$  with respect to the whole boundary. $DM$ has an universal covering space $\mathbb{H}^3$ with $\mathbb{H}^3/\Gamma=DM$, where $\Gamma\cong\pi_1(DM)$ is the covering transformation group and is discrete. Let $\pi:\H^3\rightarrow DM $ be the universal covering map.

The construction process of this section consists of three steps, and is mainly described in the unit ball model for $\mathbb{H}^3$. The first step is to study preimage of the cut locus of decorations. We will construct a $\Gamma$-invariant cell decomposition of the preimage of the cut locus and obtain some symmetry properties of the decomposition. The second step is to use the cell decomposition to construct a $\Gamma$-invariant geometric cell decomposition of $\H^3$.  And the third step is to use the universal covering and symmetry properties to induce a polyhedra decomposition of $M$.

\subsection{Cell decomposition of the preimage of the cut locus}

In this subsection, we will provide a $\Gamma$-invariant cell decomposition of preimage of the cut locus of decorations. Meanwhile, we will show that there is a  symmetry property of the cell decomposition related with geodesic boundaries, i.e. Lemma \ref{lemmaSym2}.

Let $2\nu$ be the number of cusps in $DM$. As in Section \ref{section:approach1}, we can choose a decoration $B_i$ around each cusp of $DM$, so that
\begin{itemize}
   \item[(1)] it is far away from $\partial M$;
   \item[(2)] if two cusps are symmetric along $
   \{W_j\}$, then the  decorations of these two cusps are also symmetric;
   \item[(3)] the lifting $\cup_{i=1}^{2\nu}\widetilde{B}_i$ in $\mathbb{H}^n$ are disjoint horoballs, where $\widetilde{B}_i$ is lifting of $B_i$ in $\mathbb{H}^n$.
\end{itemize}

For each pair of horoballs in $\cup_{i=1}^{2\nu}\widetilde{B}_i$, we associate a unique shortest path connecting them which is called a \emph{short cut}. Also there is an associated bisectorial geodesic plane to the short cut in $\H^3$ which is called a \emph{middle fence}. A short cut descends to the geodesic path in $DM$ from one decoration to itself or the other different decoration which is called a \emph{return path}.

\begin{definition}
The cut locus $ \textbf{C} $ of $\cup_{i=1}^{2\nu}B_i$ in $DM$ is a subset of $DM$ which consists of points that admit at least two distinct shortest paths to $\cup_{i=1}^{2\nu}B_i$.
\end{definition}

It is noted that $\cup_{j=1}^k W_j$ is belongs to $ \textbf{C} $ by the symmetry of $DM$ along $\{W_j\}$ for $1\leq j\leq k$.

Let $\widetilde{\textbf{C}} $ be the pre-image of the cut locus $ \textbf{C} $ under the universal covering map $\pi$. To construct the decomposition of $\H^3 $, we will establish a decomposition of $\widetilde{\textbf{C}} $.

A point on $ \textbf{C} $ lifts to a point on the middle fence of some short cut in $\mathbb{H}^3$. And moreover $\widetilde{\textbf{C}}$ is canonically stratified according to the number of shortest paths to the horoballs of $DM$. Then we have:

\begin{proposition}\label{pro1}
    The stratification defines a convex cellular decomposition of $\widetilde{\textbf{C}}  $.~Moreover,~the decomposition is $\Gamma$-invariant,~hence induces a cellular decomposition of $\textbf{C}  $.
\end{proposition}

\begin{proof}
Choose a component $\widetilde{U}$ of the complement of $\widetilde{\textbf{C}} $ and let $S$ be the boundary of $\overline{\widetilde{U}}$ in $\cup_{i=1}^{2\nu}\widetilde{B}_i$. Then $F=\partial\overline{\widetilde{U}} -S$ is a part of $ \widetilde{\textbf{C}} $ and formed by parts of middle fences which is called the $\emph{internal boundary}$ of $\overline{\widetilde{U}}$.

Since $DM-\cup_{i=1}^{2\nu}B_i$ is compact,~its diameter is bounded and the points on $\textbf{C}$ have bounded distance to $\cup_{i=1}^{2\nu}B_i$. The shortest arc from a point on $\textbf{C}$ to $\cup_{i=1}^{2\nu}B_i$ is lifted to geodesic arcs from some points of $\widetilde{\textbf{C}}$ to $\cup_{i=1}^{2\nu}\widetilde{B}_i$ in $\H^3$. Then the distance between $S$ and any point on $F$ is bounded. Hence the middle fences involved the points of $F$ are associated with the short cuts of bounded length.

To obtain the finiteness of return paths with bounded length, we need to study the projection of one decoration to the other with a given distance:
 \begin{claim}\label{claim1}
Suppose $\widetilde{A}$ and $\widetilde{B}$ are horoballs in different ideal points of distance $d$ in $\mathbb{H}^3$. Then the orthogonally projected image of $\widetilde{A}$ to $\widetilde{B}$ is an open metric disk of radius $\frac{1}{2}\sqrt{1+2e^{-d}}$.
 \end{claim}
\begin{proof}
Here, we use the upper half-space model of $\mathbb{H}^3$ for the proof. Without loss of generalization, let $\widetilde{A}$ and $\widetilde{B}$ be the horoballs with centers $\infty$ and $(0,0,0)$ respectively in $\mathbb{H}^3$.

Assume that the Euclidean center of $\widetilde{B}$ is $(0, 0, r)\in \mathbb{H}^3$ and the Euclidean distance from $\widetilde{B}$ to $\widetilde{A}$ is $l$, then we have $d=\log(1+\frac{l}{2r})$.

Take the geodesic plane tangent to both $\widetilde{A}$ and $\widetilde{B}$.~By directly computation,~Euclidean radius of the plane is $\sqrt{(3r+l)^2-r^2}$. Then, we can compute radius image of $\widetilde{A}$ to $\widetilde{B}$ is
an open metric disk with the radius
$$\frac{\sqrt{(3r+l)^2-r^2}}{2r+l}=\frac{1}{2} \sqrt{1+2e^{-d}}.$$
\end{proof}

By Claim \ref{claim1} and the same arguments as in the Corollary $3.3$ of Kojima \cite{Kojima}, the following holds.

\begin{corollary}
In $DM$, there exists only finitely many return paths from $\cup_{i=1}^{2\nu}B_i$ to themselves with bounded length.
\end{corollary}
Hence the middle fences involved in $F$ belong to only finitely many orbits of middle fences by the action of covering transformations preserving $S$. $F$ thus gets a locally finite invariant cellular decomposition induced
by the intersection of middle fences involved.

Using the same arguments as $\widetilde{U}$, the other components of the complement of $\widetilde{\textbf{C}}$ can also get a locally finite invariant cellular decompositions.

For two different components $\widetilde{U}_1$ and $\widetilde{U}_2$ of the complement of $\widetilde{\textbf{C}}$ such that the internal boundaries $F_1$ of $\widetilde{U}_1$ and $F_2$  of $\widetilde{U}_2$ respectively have common parts. Since the cellular decomposition of $F_i$ for $1=1,2$ is induced by the intersection of middle fences involved, the cellular decompositions of the common parts inherited from $F_1$ and $F_2$  are consistently.

By the arguments as above, $\widetilde{\textbf{C}}$ admits a cell structure which is  $\Gamma$-invariant. Therefore $\textbf{C}$ has a cellular decomposition which is induced by the cell structure of $\widetilde{\textbf{C}}$. Hence Proposition \ref{pro1} is proved.
\end{proof}

By the definition of the cut locus, every 2-cell of $\textbf{C}$  can be projected orthogonally into precisely two decorations $B_1$ and $B_2$ in $DM$ by the shortest paths respectively. In this case, we say the 2-cell \emph{faces} $B_1$ and $B_2$. Similarly, in the universal covering space $\H^{3}$, a 2-cell of $\widetilde{\textbf{C}}$ is said to \emph{face} two horoballs $B_1'$ and $B_2'$ of $\H^{3}$ if this 2-cell  can be projected orthogonally into  $B_1'$ and $B_2'$  respectively .

Then there is a symmetry property about $\textbf{C}$ that is similar to Proposition \ref{lemma2} in Section \ref{section:approach1}.

\begin{lemma}\label{lemmaSym2}

(1) If a 2-cell $f$ of $\textbf{C}$ belongs to some $W_j$ for $j=1,\ldots,k$, then $f$ faces two symmetric decorations.

(2) If a 2-cell $f$ of $\textbf{C}$ faces two decorations which are on the same side of $\cup_{j=1}^{k}W_j$, then there is another 2-cell $f'$ on the other side of $\cup_{j=1}^{k}W_j$ such that $f$ and $f'$ are symmetric along $W_j$ for $j=1,\ldots,k$.

 Moreover, the decomposition of $\textbf{C}$ is symmetric respect to $W_j~(j=1,\ldots,k)$.
\end{lemma}

\begin{proof}

Since $\widetilde{\textbf{C}}$ is constructed by parts of finitely many middle fences in $\mathbb H^3$, the dimension of $\textbf{C}$ is two. By the symmetry   of $DM$ along $\cup_{j=1}^{k}W_j$, any point of $\cup_{j=1}^{k}W_j$ admits at least two distinct shortest paths to $\cup_{i=1}^{2\nu}B_i$. Hence $\cup_{j=1}^{k}W_j$ belongs to $\textbf{C}$. Therefore, $\cup_{j=1}^{k}W_j$ consists of 2-cells (and their faces) of $\textbf{C}$.

Suppose a 2-cell $f$ of $\textbf{C}$ belongs to some $W_j$ for $j=1,\ldots,k$. Let $A, B$ be two decorations of $DM$ such that $f$ faces them. Then $A$ and $B$ are on the different sides of $\cup_{j=1}^{k}W_j$.

To derive a contradiction, we assume that $A$ and $B$ are not symmetric along any $W_j$ for $j=1,\ldots,k$. Let $p$ be a point of $f$ . We denote the shortest paths from $p$ to $A$ and from $p$ to $B$ by $l_1$ and $l_2$ respectively. By the symmetry map of $DM$ along $W_j$, there are decorations $B'$ which are the image of $B$ such that the image of $l_2$ under the symmetry is a shortest path from $p$ to $B'$. Then $p$ admits three distinct shortest paths to $\cup_{i=1}^{2\nu}B_i$, which is contradicts with the definition of $2$-cell of $\textbf{C}$. Therefore, (1) of Lemma \ref{lemmaSym2} is proved.

Suppose a 2-cell $f$ of $\textbf{C}$ faces two decorations which are on the same side of $\cup_{j=1}^{k}W_j$, then $f$ must be on the same side of $\cup_{j=1}^{k}W_j$ as the two decorations. By the symmetry map of $DM$, there is another 2-cell $f'$ on the other side of $\cup_{j=1}^{k}W_j$ such that $f$ and $f'$ are symmetric along $W_j$ for $j=1,\ldots,k$. Therefore, (2) of Lemma \ref{lemmaSym2} is proved.

Since each 2-cell of $\textbf{C}$ is either belong to $\cup_{j=1}^{k}W_j$ or is symmetric to another 2-cell along $W_j$ for $j=1,\ldots,k$, the decomposition of $\textbf{C}$ is symmetric respect to $W_j~(j=1,\ldots,k)$.
\end{proof}


\subsection{Construction of a $\Gamma$-invariant geometric decomposition of $\H^3$}

In this subsection,~we will construct a $\Gamma$-invariant geometric decomposition of $\H^3$. The decomposition is obtained by taking some duality of the cell decomposition of $\widetilde{\textbf{C}}$.  We will construct the decomposition sequentially according to their dimensions.

\paragraph*{Step 1: Construction of the edges}
We firstly construct the edges of the decomposition of $\H^3$.

Take any 2-cell of $\widetilde{\textbf{C}}$, denoted by $\mathcal{F}$.   Since $\mathcal{F}$ belongs to some middle fence according to the constructing of the cell decomposition, there are exactly two horoballs which are faced by $\mathcal{F}$. Take the shortest path of the two horoballs. Since the path is shortest, it is perpendicular to the horoballs. Then we can extend the path to a complete geodesic connecting the centers of the two horoballs. It is an edge of our decomposition.

 Do the process for the paired horoballs of all 2-cells of the cell decomposition of $\widetilde{\textbf{C}}$. Then we get the edges of the decomposition.

\paragraph*{Step 2: Construction of the faces}

We next construct the faces of the decomposition of $\H^3$.

Take any 1-cell of $\widetilde{\textbf{C}}$, denoted by $\mathcal{L}$. Since $\mathcal{L}$ is the intersection of some middle fences, one end of $\mathcal{L}$ tends to an ideal point $O\in \S^2_\infty$. There are ideal points $q_1,q_2,\ldots,q_N$ in counterclockwise order such that $q_i$ and $q_{i+1}$ is connected by a geodesic edges $l_i$ constructed as the step 1, where $i=1,2,\cdots,N$ and $q_{N+1}=q_1$. Since $\mathcal{L}$ tends to $O$, $q_1,q_2,\ldots,q_N$ are contained in a circle centered in $O$ in $\S^2_\infty$.~

Therefore,~the circle bounds a geodesic plane $H$ in the unit ball model of $\H^3$ that passes through $q_1,q_2,\ldots,q_N$ in counterclockwise order simultaneously.~

Moreover,~since $H$ passes through $q_1,q_2,\ldots,q_N$,~it intersects
orthogonally to each $B_i$ $(i=1,...,N)$ simultaneously,~which implies that short cuts between them are contained in $H$. So the edges $l_1,l_2,\cdots,l_N$ bound an ideal polygon on a geodesic plane in $\H^3 $. This is a face of our decomposition.

 Do the process for the paired horoballs of all 1-cells of the cell decomposition of $\widetilde{\textbf{C}}$. Then we get all the faces of the decomposition.
 \paragraph*{Step 3: Construction of the regions}

 We finally construct the regions of the decomposition of $\H^3$.

 Take any 0-cell of $\widetilde{\textbf{C}}$, denoted by $P$. Since $P$ is intersection of some 1-cells $\mathcal{L}_i$ of $\widetilde{\textbf{C}}$ ($i=1,2,\cdots,\mu$), there is a face $\sigma_i$ perpendicular to $\mathcal{L}_i$. We only need to show that all the faces surrounding $P$ indeed bound a region:

\begin{lemma}
     Let $\Lambda=\cup_i^\mu \overline{\sigma_i}  $ be the union of all the closed faces around $P$. Then $\Lambda$ bound a convex polyhedron in $\H^3$.
\end{lemma}
\begin{proof}
   For each $i=1,\cdots,\mu$,  $\sigma_i$ is a part of a geodesic plane $H_i$. A neighborhood of $\sigma_i$ in $\Lambda$ is contained in one side of $\H^3$ separated by $H_i$.~We call this side the \emph{inner side} of $H_i$ and call the other side the \emph{outer side} of $H_i$.

For $t\in [0,\frac{\pi}{2}]$, let $H^t_i$ be the equidistant surface in the outer side of $H_i$ at the
distance $\int_0^t \sec (\theta) d\theta $. Then for all $t\in [0,\frac{\pi}{2}]$, the ideal boundary of $H^t_i$ is the same circle in $\mathbb{S}^2_\infty$. And since $H^t_i$ is on the outer side of $H_i$ and ideal boundary is unchanged, the relationship of intersecting pairs between the $H^t_i$$(i=1,...,\mu)$ remain unchanged. So if $\sigma_i$ is bounded by $H_l$$(l=i_1,i_2,\cdots,i_k)$, then $H^t_l$$(l=i_1,i_2,\cdots,i_k)$ will bound a $\sigma^t_i$ on $H^t_i$.

Define a deformation
\begin{equation}    \begin{matrix}
    \begin{matrix}
        h_t:&\Lambda  \rightarrow \H^3\cup \S^2_\infty \\
      & \overline{\sigma_i  } \rightarrow \overline{\sigma^t_i }
        \end{matrix}
    \end{matrix} \nonumber \end{equation}

This deformation is clearly continuous and will push $\Lambda$ to the entire $\S^2_\infty$. Then $h_t(\Lambda)$ gets a stratification as $\Lambda$ for all $t$.~Note that $h_{\frac{\pi}{2}}(\Lambda)=\S^2_\infty$ bounds a convex polyhedron $\H^3$ itself.~Thus $h_0(\Lambda )=\Lambda $ bounds a convex polyhedron in $\H^3$ as well.
\end{proof}

From the above process, we construction an ideal geodesic polyhedral decomposition of $\H^3$.  Moreover, since the cell decomposition of $\widetilde{\textbf{C}}$ is $\Gamma$-invariant, we have that the geodesic decomposition constructed above is $\Gamma$-invariant.

\subsection{A geometric decomposition of $M$}

Now we have a $\Gamma$-invariant geodesic decomposition of $\H^3$. So the decomposition induces an ideal geodesic decomposition of $DM$. Now we need to show the decomposition is symmetric along $\{W_j\}$ for $1\leq j\leq k$ so that to get a mixed ideal decomposition of $M$:
\begin{proposition}
    The ideal polyhedral decomposition $DM$ obtained above will induce a mixed ideal polyhedra decomposition of $M$.
\end{proposition}

\begin{proof}
    Note that in our process of construction, the edges,faces and regions is dual to 2-cells,1-cells and 0-cells of $\widetilde{\textbf{C}}$ respectively. So after quotient by $\Gamma$, the edges,faces and regions of the decomposition of $DM$ are also dual to  2-cells,1-cells and 0-cells of $  \textbf{  C }  $ respectively.

    There are two cases to consider.

    Case $1$: If there is an ideal polyhedron $P$ such that $P$ does not intersect with any embedded submanifold $W_{j}$ for $1\leq j\leq k$ in $DM$.

    By Proposition \ref{lemmaSym2} (2), there must exist another ideal polyhedron $P'$ that is symmetric to $P$ with respect to $W_{j}$. After the symmetrical action on $DM$ along $W_j$, we obtain the original manifold $M$, and the image of $P$ and $P'$ will be two symmetric ideal polyhedra of $M$.

    Case $2$: If there is an ideal polyhedron $P''$ such that $P''$ intersects some embedded submanifold $W_{j}$ for $1\leq j\leq k$ in $DM$.

    Note that the 0-cell  dual to this polyhedral region must lie on a geodesic boundary $W_l$. For if not, take the symmetric 0-cell of it along $\{W_j\}$, then the symmetric 0-cell will be dual to the same region since $P''$ intersect some geodesic boundary. It is a contradiction since a region can be dual to only one 0-cell of $ \textbf{  C } $. Then the 1-cells around the 0-cell must be on the same geodesic boundary $W_l$ or not on geodesic boundary  since different geodesic boundary is not connected.

    It implies that $P''$  intersects only one geodesic boundary component. Moreover,  by Proposition \ref{lemmaSym2}(1), all faces intersect the geodesic boundary orthogonally. Then by the same arguments of Claim \ref{claim5}, after the symmetrical action on $DM$ along $W_j$, the quotient of $P''$ is a partially truncated polyhedron of $1$-$m$ type for some integer $m$.

\end{proof}
From the above discussion, we finally obtain a mixed ideal polyhedra decomposition of $M$ such that each
cell is either an ideal polyhedron or a partially truncated polyhedron with exactly one hyperideal vertex. Thus Theorem \ref{main} is proved in the second approach.


\section{Some comments}\label{section:comment}
\subsection{Comparison of the two proofs}\label{Section:equivalent}
The polyhedral decompositions established in the two proofs for $n=3$ have no difference. In the first proof, we project the convex hull vertically onto the one sheeted hyperboloid. So the edges in the polyhedral decomposition are projected from the bottom edges of the convex hull. Thus the edges of the polyhedral decomposition are exactly the arcs that connect the nearest horoballs, which are the same edges established in the second proof. As a consequence, the polyhedral decompositions constructed in the two proofs are exactly the same.


Although the polyhedral decompositions in the two proofs are the same, the construction processes are not completely parallel. In the first proof, we take all the arcs connecting the horoballs pairwise and project vertically to select what we actually need. But in the second proof, if we still take the convex hull construction of all the ideal points, it is difficult to directly recognize the edges that are actually needed from the convex hull construction. Therefore, we need to first construct a cell decomposition of the cut locus of the decorations. Then, we can construct the geodesic decomposition we need, which is dual to the decomposition of the cut locus.

Moreover, the symmetry properties along the totally geodesic boundaries involved in our two constructions are different. In the first proof,  we study symmetry properties of horoball centers. Meanwhile, what is important in the second proof is the symmetry property of the cell decomposition of the cut locus.

\subsection{Some furture applications}
To study the geometric and topological properties of volume finite, non-compact hyperbolic $3$-manifolds with totally geodesic boundary, a natural idea is to investigate their geometric triangulations after the mixed ideal polyhedra decomposition. By topologically triangulating the ideal polyhedra and the hyper-ideal (truncated) polyhedra, solutions to Thurston's hyperbolic gluing equations correspond to geometric triangulations. To avoid directly solving Thurston's hyperbolic gluing equations, Casson proposed an angle structure proof for triangulating $3$-manifolds. 
Based on the results in this article, in our recent work \cite{GJZ}, assuming appropriate topological conditions, we established the existence of angle structures on cusped hyperbolic $3$-manifolds with totally geodesic boundary.




\noindent
\noindent

\noindent Huabin Ge, hbge@ruc.edu.cn\\[2pt]
\emph{School of Mathematics, Renmin University of China, Beijing, 100872, P. R. China}
\\

\noindent Longsong Jia, jialongsong@stu.pku.edu.cn\\[2pt]
\emph{School of Mathematical Sciences, Peking University, Beijing, 100871, P. R. China}
\\

\noindent Faze Zhang, zhangfz201@nenu.edu.cn\\[2pt]
\emph{School of Mathematics and Statistics, Northeast Normal University, Changchun, Jilin, 130024, P.R.China} \\[2pt]

\end{document}